\documentclass[11pt]{article}
\usepackage{amsfonts}
\usepackage{mathbbold,mathrsfs,amsmath}
\usepackage{amssymb}
\usepackage{epsfig}
\usepackage{subfigure}
\usepackage{fancybox}
\usepackage{graphicx}
\usepackage{enumerate}
\newcommand{\diag}{{\rm diag}}
\newcommand{\lbar}{\overline}
\newcommand{\wdh}{\widehat}

\def\l{\left|}
\def\r{\right|}

\newcommand{\wdt}{\widetilde}
\newcommand{\e}{\varepsilon}
\newcommand{\rr}{{\Bbb R}}
\newcommand{\M}{{\cal M}}
\newcommand{\cd}{(\cdot)}
\newcommand{\nd}{\noindent}

\def\para#1{\vskip .4\baselineskip\noindent{\bf #1}}
\def\qed{\strut\hfill $\Box$}

\newtheorem{thm}{Theorem}[section]

\newtheorem{lem}[thm]{Lemma}

\newtheorem{rem}[thm]{Remark}

\newcommand{\thmref}[1]{Theorem~{\rm \ref{#1}}}
\newcommand{\lemref}[1]{Lemma~{\rm \ref{#1}}}

\def\al{\alpha}

\newcommand{\beq}[1]{\begin{equation} \label{#1}}
\newcommand{\eeq}{\end{equation}}
\newcommand{\bed}{\begin{displaymath}}
\newcommand{\eed}{\end{displaymath}}
\newcommand{\bea}{\bed\begin{array}{rl}}
\newcommand{\eea}{\end{array}\eed}
\newcommand{\ad}{&\!\!\!\disp}
\newcommand{\aad}{&\disp}
\newcommand{\barray}{\begin{array}{ll}}
\newcommand{\earray}{\end{array}}
\def\({\left(}
\def\){\right)}
\def\disp{\displaystyle}

\numberwithin{equation}{section}

\begin{document}

\title{Near-Optimal Mean-Variance Controls
under Two-time-scale Formulations and Applications\thanks{This research was
supported in part by the National Science Foundation under
CNS-1136007.}}
\author{Zhixin Yang,\thanks{Department of Mathematics,
Wayne State University, Detroit, Michigan 48202, email:
zhixin.yang@wayne.edu} \and G. Yin,\thanks{Department of
Mathematics, Wayne State University, Detroit, Michigan 48202, email:
gyin@math.wayne.edu} \and Le Yi Wang,\thanks{Department of
Electrical and Computer Engineering, Wayne State University,
Detroit, MI 48202, email: lywang@wayne.edu.} \and Hongwei
Zhang\thanks{Department of Computer Science, Wayne State University,
Detroit, MI 48202, email: hongwei@wayne.edu.}}

\maketitle

\begin{abstract}
Although the mean-variance control
was initially formulated for  financial portfolio
 management problems in which one wants to maximize
 expected return and control
 the risk, our motivations also
 stem from highway vehicle platoon
 controls that aim to maximize highway
 utility while ensuring zero accident.
This paper develops near-optimal mean-variance controls of switching
diffusion systems. To reduce the computational complexity,
with motivations from earlier work on singularly perturbed Markovian systems
\cite{SethiZ94,Yin&Zhang,GZB},
we use a two-time-scale formulation
to treat the underlying systems, which is represented by use of a
small parameter.
As the small parameter goes to 0, we obtain a
limit problem. Using the limit problem as a guide, we construct
controls for the original problem, and show that the
control so constructed is nearly optimal.

\bigskip
\nd{\bf Key Words.}
Mean-variance control, regime-switching model,
near-optimal control,
two-time scale, platoon application.

\bigskip
\nd{\bf Brief Title.}
Near-optimal Mean Variance Controls

\end{abstract}

%\newpage

\section{Introduction}
This paper
focuses on near-optimal controls of switching diffusions. Originating from
the mean-variance portfolio optimization problems, our aim
 concentrates on
reduction of computational complexity for switching diffusions where
the discrete component (switching process) has a large state space. Decomposing the
state space of the switching process into weakly connected subspaces and
aggregating the states in each subspace into
one state yield a limit system. Using the optimal controls
of the limit system, we build controls for the original systems
leading to near optimality. In addition to the traditional financial
engineering applications, our motivation stems from
 formulations of platoon
controls modeled by regime-switching systems
involving two-time-scale Markov chains, which
presents a twist of the mean-variance portfolio optimization.
This paper is written in memory of our colleague and friend
Michael Taksar, who had made significant contributions to stochastic
control, financial mathematics, and insurance risk theory; see
\cite{AT} and numerous papers Michael published.
The topic covered in the current paper is related to what Michael had
been working on. Meanwhile the application in platoon control
is a nice bifurcation from the usual
finance applications.

Although the mean-variance control
was initially formulated for  financial portfolio
 management problems in which one wants to maximize
 expected return and control
 the risk, our motivations also
 stem from highway vehicle platoon
 controls that aim to maximize highway
 utility while ensuring zero accident.
 As motivations, we identify three different but highly
 related aspects of platoon control problems
 that lead to different forms of
 the mean-variance type of problems
 that are investigated in this paper.

 First, we consider the longitudinal
 inter-vehicle distance control. To increase
 highway utility, it is desirable
 to reduce the total length of a platoon,
 which intends to reduce inter-vehicle distances.
 This strategy, however, will
 increase the risk of collision in the presence of
 vehicle traffic uncertainties.
 A tradeoff of these factors leads to a desired nominal length.
 Deviation from this nominal platoon length compromises
 either safety or highway utility.
 Since vehicle movements are subject to many
 random factors on road, weather, and traffic conditions,
 the total platoon length is actually a stochastic process.
 The desire to control the platoon length
 toward its designated target (its mean)
 with small deviations (its variance)
 can be mathematically modeled as a mean-variance
 optimization problem with subsystem states as inter-vehicle distances.

 The second scenario is communication resource
 allocation of bandwidths for vehicle to vehicle (V2V)
 communications. For a given maximum throughput
 of a platoon communication system,
 the communication system operator
 must find a way to assign this resource
 to different V2V channels.
 If the total bandwidth used is
 lower than assigned bandwidth,
 there will be resource waste.
 Conversely, usage of bandwidths over
 the budget may incur high costs or interfere with
 other platoons' operation. In this case,
 each channel's bandwidth usage is the state
 of the subsystem. Their summation is a random process
 and is desired to approach the maximum throughput
 (the desired mean at the terminal time)
 with small variations.
 Consequently, it becomes a mean-variance control problem.

 Finally, we may view platoon fuel consumption
 (or similarly, total emission).
 The platoon fuel consumption is
 the summation of vehicle fuel consumptions.
 Due to variations in vehicle sizes and speeds,
 each vehicle's fuel consumption is a controlled
 random process.
 Tradeoff between a platoon's team acceleration/maneuver
 capability and fuel consumption can be summarized in
 a desired platoon fuel consumption rate.
 Assigning fuels to different vehicles result
 in coordination of vehicle operations modeled
 by subsystem fuel rate dynamics.
 To control the platoon fuel consumption
 rate to be close to the designated value,
 one may formulate this as a mean-variance control problem.

Although the MV approach has never been applied to platoon
 controls, it has distinct advantages:
 1) unlike heuristic methods such as neural network
 optimization and genetic algorithms,
  the MV method is simple but rigorous; 2)
  the MV method is computationally efficient;
   3) the form of the solution (i.e., efficient frontier)
   is readily applicable
   to assessing risks in platoon formation, hence is practically appealing.

The origin of the mean-variance optimization problem can be traced
back to the Nobel-price-winning work of Markowitz \cite{Mark52}.
The salient feature of the model is that, in the context of
finance, it enables an investor to seek desired expected return after
specifying the acceptable risk level quantified by the
variance of the return. The mean-variance approach has become the
foundation of modern finance theory and has inspired numerous
extensions and applications. Using the ideas of backward stochastic
differential equations, the mean-variance
problem for a continuous-time model was studied in \cite{ZhouL}. Note that
in the mean variance problems,
the matrix related to the control (known as control weight)
 is not positive definite.
To take into consideration of random environments not representable
using the usual stochastic differential equation setup, we developed more
complex
 models with random switching in \cite{Y&Z}.

In this paper, we consider the case that the
random process representing discrete events (the environment)
has a large state space.
The physical system is such that not all of the discrete event states
change at the same rate.
Some of them vary rapidly and others change slowly. The fast and
slow variations are in high contrast
resulting in a two-time-scale formulation.
Taking advantage of the time-scale separation,
we use an averaging
approach to analyze the system, which
largely explores the weak and strong interactions
of the switching diffusion due to
the Markov chain.
 The rationale is to
aggregate the states according to their jump rates and replace the
actual system with its average. Using optimal control of the limit
problem as a bridge, we then construct controls for the original
systems leading to feasible
approximation schemes.
In \cite{LiuYZ}, we treated a class of LQ problems with switching
 by concentrating on
the associated Riccati systems of equations, whereas in this paper,
we focus on mean-variance controls and examining certain associated
systems of differential equations.
We consider the case that the
Markov chains have recurrent states as well as inclusion of transient states. 
These
approximation schemes give us nearly optimal controls. Focusing
on approximated optimality, we succeed in reducing the complexity
of the underlying systems substantially.

The rest of paper is arranged as follows. Section \ref{sec:form}
begins with the
formulation of the two-time-scale platoon problems. Section
\ref{sec:prelim} proceeds with the study of the
underlying mean-variance problem. Using
completing square techniques, we derive the corresponding Riccati
equations and optimal control  for the non-definite control problem.
Section \ref{sec:near}
focuses on near-optimal controls of the mean-variance problems.
First, Markov
chains with recurrent states are treated
and then  inclusion of
transient states are considered. Using probabilistic arguments and analytic
techniques, the approximation schemes are shown to be nearly
optimal.  Finally, we
conclude this paper with further thoughts and additional remarks in
Section \ref{sec:rem}.

\section{Formulation}\label{sec:form}
 We begin with a complete probability space $(\Omega,
\mathcal{F}, P)$. Consider
a time-homogeneous Markov
chain in continuous time
 taking values in the state space
 $\mathcal{M}=\{1,2,\ldots,m\}$ and a standard
$d$-dimensional standard Brownian motion
$w(t)=(w_1(t),w_2(t),\ldots,w_d(t))'$ (where $a'$ denotes the
transpose of $a \in \rr^{l_1\times l_2}$ with $l_i \ge 1$) that is
independent of the Markov chain $\al(t)$. The Markov chain is used
to represent discrete events and random environment etc. In
\cite{Y&Z}, we considered the Markovitz's mean-variance portfolio
selection problem in which the environment is allowed to vary
randomly leading to a regime-switching model. In this paper, we
continue using the setup as in that of \cite{Y&Z}. In addition to
the finance applications, we have in mind the platoon control
problems as mentioned in the introduction. Mathematically, the new
feature considered here is that the state space of the discrete
event process $\al\cd$ is large. Treating mean-variance control
problems thus requires handling of large-scale systems. Such a case
naturally arises in the networked system formulation. The large
scale feature however renders the optimal control a difficult task.
To reduce the computational complexity, we note that for the discrete event
process (the Markov
chain), not all states are varying at the same rate. Some clusters of
states vary rapidly and others change slowly. Using the relative
transition rates, we decompose the state space $\M$ into subspaces
$\M_i$ such that within each $\M_i$, the
transitions happen frequently and among different clusters the
transitions are relatively infrequently. To reflect the different
transition rates, we let $\al(t)=\al^\e(t)$ where $\e>0$ is a small
parameter so that the generator of the Markov chain is given by
\beq{mar}%
Q^\e=\frac{ 1}{\e}\wdt Q+\wdh{Q}. \eeq Define
$\mathcal{F}_t=\sigma\{W(s),\al^\e(s):0\le s\le t\}.$
Denote
\beq{mar1}%
 Q^\e f(\cdot)(i)=\sum_{j\neq
i}q^\e_{ij}(f(\cdot,j)-f(\cdot,i))\eeq for a suitable $f(\cdot)$.
Suppose that $x_i^\e\cd$ are real-valued functions with $i=0,\ldots,
d_1$
such that
\beq{0}\barray dx_0^\e(t)\ad=r(t,\al^\e(t))x_0^\e(t)dt\\
x_0^\e(0)\ad=x_0, \ \al^\e(0)=\al\earray\eeq for
$\al^\e(t)\in \{1,2,\ldots, m\}$. The flows of the other $d_1$ nodes follow
geometric Brownian motion: \beq{2}\barray
dx_i^\e(t)\ad=x_i^\e(t)r_i(t,\al^\e(t))dt+x_i^\e(t)
\sigma_i(t,\al^\e(t))dw(t)\\
x_i^\e(0)\ad=x_i,\ \al^\e(0)=\al \text{ for }
i=1,2,\ldots,d_1, \ \al \in \M,\earray \eeq where $\sigma_i(t,\al^\e(t))=
(\sigma_{i1}(t,\al^\e(t)),
\sigma_{i2}(t,\al^\e(t)),\ldots,\sigma_{id}(t,\al^\e(t)))\in
\rr^{1\times d}$. In the finance application, $x^\e_0\cd$ represents
an investor's bank account value, whereas $x^\e_i\cd$ for each
$i=1,\ldots, d_1$ is his wealth devoted to the $i$th
stock. In the networked control problems, we use $x^\e_i\cd$ to
represent the flows of the $i$th node. We can represent the wealth
of the investor or the total flows of the entire system as
 $x^\e(t)$.

 For consistency with the current literature on the
 MV problems, we shall still use the term
 ``portfolio" in our network problems.
 In the traditional market analysis setting,
a portfolio is a vector consisting of the dollar
values of different stocks.
When applied to our network systems, a portfolio will be
 the vector of inter-vehicle distances in platoon control,
 or individual channel throughput in communication resource allocation, or
 individual vehicle fuel consumption in platoon fuel management.
The portfolio
 selection involves finding the strategy to
select the proportion $n_i\cd$ of the $i$th stock investment.
Similarly,
for the platoon problem,
 we need to decide the
proportion $n_i(t)$ of the flow $x_i^\e(t)$ on node $i$.
In these cases, we denote their sum as
 \bea \disp x^\e(t)=\sum^{d_1}_{i=0}n_i(t)x_i^\e(t).\eea
 Then we have
\beq{x}%
\barray%
dx^\e(t)
\ad=[r(t,\al^\e(t))x^\e(t)+B(t,\al^\e(t))u(t)]dt
+u'(t)\sigma(t,\al^\e(t))dw(t), \ t\in [0,T],\\
x^\e(0)\ad=\hat{x}=\sum^{d_1}_{i=1}n_i(0)x_i ,\ \
\al^\e(0)=\al,\earray \eeq where
$$B(t,\al^\e(t))=(r_1(t,\al^\e(t))
-r(t,\al^\e(t)),r_2(t,\al^\e(t))-r(t,\al^\e(t)),
\ldots,r_{d_1}(t,\al^\e(t))-r(t,\al^\e(t))),$$ $\sigma(t,\al^\e(t))= (
  \sigma _1 (t,\al^\e(t)),\ldots, \sigma _{d_1}
  (t,\al^\e(t)))'
\in\rr^{d_1\times 1},$ $u(t)=(u_1(t),\dots,
u_{d_1}(t))'\in\rr^{d_1\times 1}$, and $u_i(t)=n_i(t)x_i(t)$ where
$u_i(t)$ is the total amount of in the $i$th stock or the amount of
flow for node $i$ at time $t$, for $i=1,2,\ldots,d_1$. We assume
throughout this paper that $r(t,i)$, $B(t,i)$, $\sigma(t,i)$ are
measurable and uniformly bounded in $t$ and we also assume the
non-degeneracy condition is satisfied, i.e., there is some $\delta
>0$, $a(t,i)=\sigma(t, i)\sigma'(t,i)\ge \delta I$
for any $t\in[0,T]$ and each $i\in \M$. We denote by
$L^2_{\mathcal{F}}(0,T;\mathbb{R}^{l_0})$ the set of all
$\mathbb{R}^{l_0}$ valued, measurable stochastic processes $f(t)$
adapted to $\{\mathcal{F}_t\}_{t\ge 0}$ such that
$E\int^T_0|f(t)|^2dt< \infty$.

Let $\mathcal{U}$, the set of control,  be a compact subset in
$\mathbb{R}^{d_1\times 1}$. The $u(\cdot)$ is said to be admissible
if for a $\mathcal{U}$ valued control $u(\cdot)\in
L^2_{\mathcal{F}}(0,T;\mathbb{R}^{d_1})$, the equation \eqref{x} has
a unique solution $x^\e(\cdot)$ corresponding to $u(\cdot)$. In this
case, we refer to $(x^\e(\cdot),u(\cdot))$ as an admissible pair.
 Our objective is to find an admissible
 control $u(\cdot)$ among all the admissible controls
given that expected terminal flow value of the whole system is
$Ex^\e(T)=z$ for some given $z\in \mathbb{R}$, so that the risk
measured by the variance of terminal of the flow is minimized.
Thus, we have the following goal \beq{obj}
\barray \ad
\hbox{min }
J(x,\al,u(\cdot))=E[x^\e(T)-z]^2\\
\ad \hbox{subject to } Ex^\e(T)=z. \earray \eeq
Recall that the problem is called feasible if there is at
least one portfolio satisfying
all the constraints. The problem is called finite if it is feasible
and the infimum of
$J(x, \al, u(\cdot))$ is finite.
An optimal portfolio to the above problem, if it
ever exists, is called an efficient portfolio corresponding to $z$, and
the corresponding
$(Var x(T), z)\in \rr^2$  is called an
efficient point,
 The set of all
the efficient points
is called the efficient frontier.
The solution of the problem can be
obtained by using the result of \cite{Y&Z}. In fact, we can obtain
the efficient frontier as well as the so-called mutual fund
theorems. This however is not the end of the story but rather the
starting point of the current paper. In this paper, we consider the
case that $|\M|=m$ is large, thus we have to solve a system of $m$
equations where $m$ is large. Computationally, this is rather
cumbersome. Therefore, our effort is devoted to reducing the
complexity.

\section{Preliminary Results}\label{sec:prelim}
This section presents preliminary results concerning the solutions of
systems. The results include feasibility, existence and uniqueness of
the solution, and continuity.
For the feasibility part of our problem, we present the following
lemma. The detailed proof can be found in \cite[Theorem3.3]{Y&Z}.
\begin{lem}
The mean variance problem \eqref{obj} is feasible for every $z\in
\rr$  if and only if \beq{fe}\barray
E\int^T_0|B(t,\al^\e(t))|^2dt>0.\earray\eeq
\end{lem}
 Now let us proceed to the
study of optimality. To handle the constraint part in problem
\eqref{obj}, we apply the Lagrange multiplier technique and get
unconstrained problem (see, e.g.,\cite{Y&Z}) with multiplier
$\lambda\in \mathbb{R}$:
\beq{obj2}%
\barray%
\ad \hbox{min }
J(x,\al,u(\cdot),\lambda)=E[x^\e(T)+\lambda-z]^2-\lambda^2\\
\ad \hbox{subject to } E x^\e(T)=z, \
\hbox{ with }\ (x^\e(\cdot),u(\cdot)) \text{ admissible. }
\earray \eeq To find the minimum of $J(x,\al, u\cd, \lambda)$, it
suffices to choose $u\cd$ so that $E(x^\e(T)+\lambda-z)^2$ is
minimized. We regard this part as $J^\e(x,\al,u\cd)$ in what
follows. In this section, we will proceed to solve the unconstrained
problem \eqref{obj2}. Let
$v^\e(x,\al)=\inf_{u(\cdot)}J^\e(x,\al,u(\cdot))$ be the value
function. First define
\beq{r}%
\barray \rho(t,i)=B(t,i)[\sigma(t,i)\sigma'(t,i)]^{-1}B'(t,i),\ i\in
\{1,2,\ldots,m\}. \earray \eeq
 Consider the following
two systems of ODEs for $i=1,2,\ldots, m$: \beq{p}\barray
\ad\dot{P}^\e(t,i)=P^\e(t,i)[\rho(t,i)-2r(t,i)]
-\sum^m_{j=1}q^\e_{ij}P^\e(t,j)\\
\ad P^\e(T,i)=1 .\earray \eeq and
 \beq{h}\barray
\dot{H}^\e(t,i)=\ad H^\e(t,i)r(t,i)
-\frac{1}{P^\e(t,i)}\sum^m_{j=1}q^\e_{ij}P^\e(t,j)H^\e(t,j)
%\\ \ad
+\frac{H^\e(t,i)}{P^\e(t,i)}\sum^m_{j=1}q^\e_{ij}P^\e(t,j)\\
H^\e(T,i)=\ad 1. \earray \eeq

The existence and uniqueness of solutions to the above two systems
of equations are evident as both are linear with uniformly bounded
coefficients. Applying the generalized It\^{o}'s formula to
\bea%
v^\e(t,x^\e(t),i)=P^\e(t,i)(x^\e(t)+(\lambda-z)H^\e(t,i))^2 \eea
 and using the completing square techniques, we obtain
\beq{eq1}\barray
\ad\! d P^\e(t,i)[x^\e(t)+(\lambda-z)H^\e(t,i)]^2\\
\aad
=2P^\e(t,i)[x^\e(t)+(\lambda-z)H^\e(t,i)]dx^\e(t)+P^\e(t,i)(dx^\e(t))^2
\\
\aad \ +\sum^m_{j=1}q^\e_{ij}P^\e(t,j)
[x^\e(t)+(\lambda-z)H^\e(t,j)]^2dt\\
\aad \
+\dot{P}^\e(t,i)[x^\e(t)+(\lambda-z)H^\e(t,i)]^2dt
+2P^\e(t,i)[x^\e(t)+(\lambda-z)H^\e(t,i)]
(\lambda-z)\dot{H}^\e(t,i)dt.\earray\eeq
Therefore, by plugging in the dynamic equation satisfied by $
P^\e(t,i)$ and $H^\e(t,i)$, we have the following expression:
\beq{v}%
\barray%
\ad
dP^\e(t,i)[x^\e(t)+(\lambda-z)H^\e(t,i)]^2\\
\ad=P^\e(t,i)\{u'(t)\sigma(t,i)\sigma'(t,i)u(t)+2u'(t)B'(t,i)
[x^\e(t)+(\lambda-z)H^\e(t,i)]\\
\aad+2r(t,i)x^\e(t)[x^\e(t)+(\lambda-z)H^\e(t,i)]\}dt
-\sum^m_{j=1}q^\e_{ij}P^\e(t,j)[x^\e(t)+(\lambda-z)H^\e(t,i)]^2dt\\
\aad+2P^\e(t,i)[x^\e(t)+(\lambda-z)H^\e(t,i)] (\lambda-z)\{
H^\e(t,i)r(t,i)-\frac{1}{P^\e(t,i)}\sum^m_{j=1}
q^\e_{ij}P^\e(t,j)H^\e(t,j)\\
\aad+\frac{H^\e(t,i)}{P^\e(t,i)}\sum^m_{j=1}
q^\e_{ij}P^\e(t,j)\}dt+\sum^m_{j=1}q^\e_{ij}P^\e(t,j)
[x^\e(t)+(\lambda-z)H^\e(t,j)]^2dt\\
\aad+[\rho(t,i)-2r(t,i)]P^\e(t,i)[x^\e(t)
+(\lambda-z)H^\e(t,i)]^2dt+(\cdots)dw(t)\\
\ad=
P^\e(t,i)\{(u(t)+(\sigma(t,i)\sigma'(t,i))^{-1}
B'(t,i)[x^\e(t)+(\lambda-z)H^\e(t,i)])'[\sigma(t,i)\sigma'(t,i)]\\
\aad \ \times
(u(t)+(\sigma(t,i)\sigma'(t,i))^{-1}
B'(t,i)[x^\e(t)+(\lambda-z)H^\e(t,i)])\}dt\\
\aad \
+(\lambda-z)^2\sum^m_{j=1}q^\e_{ij}P^\e(t,j)
[H^\e(t,j)-H^\e(t,i)]^2dt+(\cdots)dw(t).\earray
\eeq Integrating both sides of the above equation from $0$ to $T$
and taking expectation,
we obtain%
\beq{ex}%
\barray%
\ad E[x^\e(T)+\lambda-z]^2\\
\aad \ = P^\e(0,\al)[x+(\lambda-z)H^\e(0,\al)]^2\\
\aad\quad +E\int^T_0(\lambda-z)^2\sum^m_{j=1}
q^\e_{ij}P^\e(t,j)[H^\e(t,j)-H^\e(t,i)]^2dt\\
\aad\quad
+E\int^T_0P^\e(t,i)[u(t)-u^{\e,*}(t)]'[\sigma(t,i)\sigma'(t,i)]
[u(t)-u^{\e,*}(t)]dt.
\earray%
\eeq%

Thus, the optimal control $u^{\e,*}$ has the form
\beq{op}
u^{\e,*}(t,\al^\e(t),x^\e(t))=-(\sigma(t,\al^\e(t))\sigma'(t,\al^\e(t)))^{-1}
B'(t,\al^\e(t))[x^\e(t)+(\lambda-z)H^\e(t,\al^\e(t))].\eeq
 Now We
introduce the following two lemmas here for the subsequent use.

\begin{lem}\label{bdd}
  The solution of equations \eqref{p} and \eqref{h} satisfy
  $0<P^\e(t,i)\le c$ and $0<H^\e(t,i)\le 1$ for all $t\in
  [0,T],i=1,2,\ldots,m$.
\end{lem}

\para{Proof.}
For the $H^\e(t,i)$, by employing the idea similar to
\cite[Proposition 4.1]{Y&Z}, we can get the claim above. Here, we
consider the case for $P^\e(t,i)$. First, by applying a variation of
constant formula to \eqref{p} we have
\bea \disp
P^\e(t,i)=e^{-\int^T_t[\rho(s,i)-2r(s,i)-q^\e_{ii}]ds}+\int^T_t
e^{-\int^s_t[\rho(\tau,i)-2r(\tau,i)-q^\e_{ii}]d\tau}\sum^m_{j\neq
i}q^\e_{ij}P^\e(s,j)ds.\eea Construct a Picard sequence of
$P^\e_k(\cdot,i)$ for $t\in[0,T], i=1,2,\ldots,m,k=0,1,\ldots$ as
follows
\bea P^\e_0(t,i)\ad=1,\\
P^\e_{k+1}(t,i)\ad=e^{-\int^T_t[\rho(s,i)-2r(s,i)-q^\e_{ii}]ds}+\int^T_t
e^{-\int^s_t[\rho(\tau,i)-2r(\tau,i)-q^\e_{ii}]d\tau}\sum^m_{j\neq
i}q^\e_{ij}P^\e_k(s,j)ds\eea

Noting that $q^\e_{ij}\ge 0$ for all $j \neq i$, we have for
$k=0,1,\ldots$ \bea P^\e_k(t,i)\ge
e^{-\int^T_t[\rho(s,i)-2r(s,i)-q^\e_{ii}]ds}>0, \eea Realizing that
$P^\e(t,i)$ is the limit of the Picard sequence $P^\e_k(t,i)$ as $k
\to \infty$. Thus, $P^\e(t,i) >0$. To get the upper bound, we first
consider the bounds of value function $v^\e(x,\al)$. Clearly,
$v^\e(x,\al)\ge 0$ since $J^\e(x,\al,u(\cdot))\ge 0$ for all
admissible $u(\cdot)$. We choose $u_0(t)=-ax^\e(t)$, $a$ is a
nonzero vector in $\mathbb{R}^{d_1}$ and $x^\e(t)=\tilde{x}$, then
we have $E(x^\e(T))^2\le \tilde{x}^2+k\int^T_t E(x^\e(s))^2ds $
according to It\^{o}'s formula. We further have $E(x^\e(T))^2\le
\tilde{x}^2e^{kT}$ by virtue of Gronwall's inequality for all
$t\in[0,T]$. Now, note that for $0\le t\le T$, \bea
v^\e(\tilde{x},i)\le J^\e(\tilde{x},i,u(\cdot))\le
E[x^\e(T)+\lambda-z]^2\le 2\tilde{x}^2e^{kT}+2(\lambda-z)^2.\eea
Then we have \bea P^\e(t,i)(\tilde{x}+(\lambda-z)H^\e(t,i))^2\le
2\tilde{x}^2e^{kT}+2(\lambda-z)^2. \eea Dividing both
 sides of this inequality by $\tilde{x}^2$ and setting
$\tilde{x}\to\infty$, we have $P^\e(t,i)\le 2e^{kT}$. \qed

\begin{lem}\label{unc}
For $i\in\mathcal{M}$, the solutions of \eqref{p} and \eqref{h} are
uniformly Lipschitz on $[0,T]$.
\end{lem}

\para{Proof.}
Let us just consider the part of $P^\e(t,i)$ since the proof for the
case of $H^\e(t,i)$ is similar. Given that the solution for equation
\eqref{p} is
 \bea P^\e(t,i)\ad =
e^{-\int^T_t[\rho(s,i)-2r(s,i)]ds}+\int^T_t
e^{-\int^s_t[\rho(\tau,i)-2r(\tau,i)]d\tau}\sum^m_{j=1}q^\e_{ij}P^\e(s,j)ds\\
\ad =
e^{-\int^{t+\Delta}_t[\rho(s,i)-2r(s,i)]ds
-\int^T_{t+\Delta}[\rho(s,i)-2r(s,i)]ds}\\
\aad \ +\int^{t+\Delta}_t
e^{-\int^s_t[\rho(\tau,i)-2r(\tau,i)]d\tau}\sum^m_{j=1}q^\e_{ij}P^\e(s,j)ds\\
\aad \ +\int^T_{t+\Delta}
e^{-\int^{t+\Delta}_t[\rho(\tau,i)-2r(\tau,i)]d\tau-\int^s_{t+\Delta}[\rho(\tau,i)-2r(\tau,i)]d\tau}
\sum^m_{j=1}q^\e_{ij}P^\e(s,j)ds \eea Given that \bea
P^\e(t+\Delta,i)=\ad
e^{-\int^T_{t+\Delta}[\rho(s,i)-2r(s,i)]ds}+\int^T_{t+\Delta}
e^{-\int^s_{t+\Delta}[\rho(\tau,i)-2r(\tau,i)]d\tau}\sum^m_{j=1}q^\e_{ij}P^\e(s,j)ds.\eea
Then we have \bea \ad
|P^\e(t,i)-P^\e(t+\Delta,i)|\\
\aad \ = |e^{-\int^T_{t+\Delta}[\rho(s,i)-2r(s,i)]ds}(1-e^{-\int^{t+\Delta}_t[\rho(s,i)-2r(s,i)]ds})\\
\aad\quad + \int^T_{t+\Delta}e^{-\int^s_{t+\Delta}[\rho(\tau,i)-2r(\tau,i)]d\tau}(1-e^{-\int^{t+\Delta}_t[\rho(\tau,i)-2r(\tau,i)]
d\tau})\sum^m_{j=1}q^\e_{ij}P^\e(s,j)ds\\
\aad \quad +\int^{t+\Delta}_t
e^{-\int^s_t[\rho(\tau,i)-2r(\tau,i)]d\tau}\sum^m_{j=1}q^\e_{ij}P^\e(s,j)ds|.\eea
As $\Delta \to 0$, $$
1-e^{-\int^{t+\Delta}_t[\rho(\tau,i)-2r(\tau,i)] d\tau} \to 0$$
and $$ \int^{t+\Delta}_t
e^{-\int^s_t[\rho(\tau,i)-2r(\tau,i)]d\tau}\sum^m_{j=1}q^\e_{ij}P^\e(s,j)ds
\to 0$$  hold for any $t \in [0,T]$ uniformly, therefore,
$P^\e(t,i)$ is uniformly Lipschitz on $[0,T]$. \qed

Due to the large dimensionality, it is highly computation intensive
to obtain
the optimal controls.
To overcome the difficulty, we device a near-optimal control scheme.
We will show that as $\e\to 0$, there is a limit problem. For the
limit problem, we can obtain optimal controls as given in
\cite{Y&Z}. Then we use the optimal control of the limit problem
to construct controls of the original problem and show that the
constructed control is asymptotically optimal.

\section{Near-Optimal Controls}\label{sec:near}
\subsection{Recurrent States}
Assume that $\wdt{Q}$ can be put into a block-diagonal form
$\wdt{Q}=\diag(\wdt{Q}^1,\ldots, \wdt{Q}^l)$ in which $\wdt{Q}^k\in
\mathbb{R}^{m_k \times m_k}$ are irreducible for $k=1,2,\ldots,l$
and $\sum^l_{k=1}m_k=m $. $\wdt{Q}^k$ denotes the $kth$ block matrix
in $\wdt{Q}$. Denote by $\mathcal{M}_k=\{s_{k1},s_{k2},\ldots,s_{km_k}\}$
 the states corresponding to $\wdt{Q}^k$ and note
$$\mathcal{M}=\mathcal{M}_1\cup\mathcal{M}_2 \cup\cdots
\cup\mathcal{M}_l.$$
Note that the $\wdt{Q}^k=(\tilde{q}^k_{ij})_{m_k \times m_k}$ and
$\wdh{Q}=(\wdh{q}_{ij})_{m\times m}$ are generators.

The slow and fast components are coupled through weak and strong
interactions in the sense that the underlying Markov chain
fluctuates rapidly within a single group $\mathcal{M}_k$ and jumps
less frequently among groups $\mathcal{M}_k$ and $\mathcal{M}_j$ for
$k\neq j$.
By aggregating the states in $\mathcal{M}_k$ as one state $k$, we
obtain an aggregated process $\lbar{\al}^\e(\cdot)$ defined by
$\lbar{\al}^\e(t)=k$ when $\al^\e(t)\in \mathcal{M}_k$.
Although $\lbar \al^\e(t)$ is generally not Markovian, by virtue of
\cite[Theorem7.4]{Yin&Zhang}, $\lbar{\al}^\epsilon(\cdot)$
converges weakly to a Markov chain $\lbar{\al}(\cdot)$ with generator
$\overline{Q}=(\lbar{q}_{ij})$ satisfying
\bea%
\overline{Q}=\diag(\mu^1,\mu^2,\ldots,\mu^l)\wdh{Q}\diag(\mathbbold{1}_{m_1},
\mathbbold{1}_{m_2},\ldots, \mathbbold{1}_{m_l}).\eea where $\mu^k$
is the stationary distribution associated with
$\wdt{Q}^k,k=1,2,\ldots,l$, and $\mathbbold{1}_n=(1,1,\ldots,1)\in
\mathbb{R}^{n\times 1}$. For subsequent use, we define
$\lbar{F}(t,k)=\sum^{m_k}_{j=1}\mu^k_jF(t,s_{kj})$ for
$F(t,s_{kj})=r(t,s_{kj}), B(t,s_{kj})$ and $\rho(t,s_{kj})$. The
following theorems are concerned with the convergence and nearly
optimal control.

\begin{thm}\label{4.1}
For $k=1,2,\ldots,l$ and $j=1,2,\ldots,m_k$, $P^\e(t,s_{kj})\to
\overline{P}(t,k)$ and $H^\e(t,s_{kj})\to \overline{H}(t,k)$
uniformly on $[0,T]$ as $\e \to 0$, where $\overline{P}(t,k)$ and
$\overline{H}(t,k)$ are the unique solutions of the following
differential equations for $k=1,2,\ldots,l$, \beq{P2} \barray
\dot{\overline{P}}(t,k)=\ad(\overline{\rho}(t,k)-2\lbar{r}(t,k))\overline{P}(t,k)-
\bar{Q}\overline{P}(t,\cdot)(k)\\
\overline{P}(T,k)=\ad 1. \earray \eeq and \beq{H2}\barray
\dot{\overline{H}}(t,k)=\ad\overline{r}(t,k)\overline{H}(t,k)-\frac{1}{\overline{P}(t,k)}
\bar{Q}\overline{P}(t,\cdot)\overline{H}(t,\cdot)(k)
%\\ \ad
+\frac{\overline{H}(t,k)}{\overline{P}(t,k)}\bar{Q}\overline{P}(t,\cdot)(k)\\
\overline{H}(T,k)=\ad 1. \earray\eeq
\end{thm}

\para{Proof.}
We prove the convergence of $P^\e$  (the proof of $H^\e$ is similar). By
virtue of \lemref{bdd} and \lemref{unc}, $P^\e(t,s_{kj})$
is equicontinuous and uniformly bounded, it follows from
Arzela-Ascoli theorem that, for each sequence of $\e \to 0$, a
further subsequence exists (we still use the index $\e$ for the sake
of simplicity) such that $P^\e(t,s_{kj})$ converges uniformly on
$[0,T]$ to a continuous function, say, $P^0(t,s_{kj})$. First, we
show $P^0(t,s_{kj})$ is independent of $j$. Given that
$$
P^\e(t,s_{kj})=1-\int^T_t[P^\e(s,s_{kj})(\rho(s,s_{kj})
-2r(s,s_{kj}))-Q^\e P^\e(s,\cdot)(s_{kj})]ds .$$ Multiplying both
sides of above equation by $\e$ yields that \bea \disp 0=\lim_{\e \to
0}\int^T_t\wdt{Q}^kP^\e(s,\cdot)(s_{kj})ds=\int^T_t\wdt{Q}^kP^0(s,\cdot)(s_{kj})ds.\eea
Thus, in view of the continuity of $P^0(t,\cdot)(s_{kj})$, we obtain
\beq{ct}\barray \wdt{Q}^kP^0(t,\cdot)(s_{kj})=0
 \text{ for }t\in [0,T].\earray\eeq

  Given the fact that $\wdt{Q}^k$ is
irreducible, we have $P^0(t,s_{kj})=P^0(t,k)$ which is independent
of $j$. Now let us multiply $P^\e(t,s_{kj})$ by $\mu^k_j$ and then
add the index $j$. Recall the definition of $\lbar{F}(t,k)$, we have
the following equation \bea \disp
\sum^{m_{k}}_{j=1}\mu^k_jP^\e(t,s_{kj})
\ad=1-\sum^{m_{k}}_{j=1}\mu^k_j\int^T_t[P^\e(s,s_{kj})(\rho(s,s_{kj})
-2r(s,s_{kj}))-Q^\e P^\e(s,\cdot)(s_{kj})]ds
. \eea Letting
$\e\to 0$ and noting that uniform convergence of $P^\e(t,s_{kj})\to
P^0(t,k)$ and $\mu^k$ is the stationary distribution corresponding
to $\wdt{Q}^k$, we have \bea \disp
(\sum^{m_k}_{j=1}\mu^k_j\wdh{Q}\mathbbold{1}_{m_k})P^0(t,\cdot)(k)=\lbar{Q}P^0(t,\cdot)(k).\eea
Therefore, we obtain%
\bea%
P^0(t,k)=\ad 1-\int^T_t\Big(P^0(s,k)(\lbar{\rho}(s,k)
-2\lbar{r}(s,k)-\lbar{Q}
P^0(s,\cdot)(k)\Big)ds%
\eea%
Then the uniqueness of solution of the Riccati equation implies
$P^0(s,k)=\lbar{P}(s,k)$. Therefore, $P^\e(t,s_{kj})\to \lbar{P}(t,k)$
and the proof is thus concluded. \qed

It follows that $P^\e(t,s_{kj})\to \lbar{P}(t,k)$ and
$H^\e(t,s_{kj})\to \lbar{H}(t,k)$. We thus have $v^\e(t,s_{kj},x)\to
\lbar{v}(t,k,x)$ as $\e \to 0$, in which
$\lbar{v}(t,k,x)=\lbar{P}(t,k)(x+(\lambda-z)\lbar{H}(t,k))^2$, where
$\lbar{v}(t,k,x)$ corresponds to the value function of a limit
problem. Let $\mathcal{U}$ denote the control set for the limit
problem:
$\mathcal{U}=\{U=(U^1,U^2,\ldots,U^l):U^k=(u^{k1},u^{k2},\ldots,u^{km_k}),u^{kj}\in
\mathbb{R}^{d_1}\}$. Define
\bea \ad f(t,x,k,U)=\sum^{m_k}_{j=1}\mu^k_jr(t,s_{kj})x
+\sum^{m_k}_{j=1}\mu^k_jB(t,s_{kj})u^{kj}(t) \ \hbox{
and }\\
\ad g(t,k,U)=(g_1(t,k,U),\ldots, g_d(t,k,U)) \text{ with }
g_i(t,k,U)=\sqrt{\sum^{m_k}_{j=1}{\mu^k_j}(\sum^{d_1}_{n=1}u^{kj}_n\sigma_{ni}(t,\al^\e(t)))^2}.\eea
Recall that $\sigma(t,\al^\e(t))=(\sigma_{ni}(t,\al^\e(t)))\in
\rr^{d_1\times d}$ and note that $u^{kj}_n$ is the $n$th component
of the $d_1$-dimensional variable. The corresponding dynamic system
of the state is \beq{li}  dx(t)=f(t,x(t),\lbar{\al}(t),
U(t))dt+\sum^d_{i=1}g_i(t,\lbar{\al}(t),U(t))dw_i(t).\eeq where
$\overline{\al}(\cdot)\in \{1,2,\ldots,l\}$ is a Markov chain
generated by $\lbar{Q}$ with $\lbar{\al}(0)=\al$. Calculation
similar to \eqref{eq1} and \eqref{v} shows that the optimal control
for this limit problem is \bea
\ad U^*(t)=(U^{1*}(t,x),U^{2*}(t,x),\ldots,U^{l*}(t,x)), \ \hbox{ with }
 U^{k*}(t,x)=(u^{k1*}(t,x),u^{k2*}(t,x),\ldots,u^{km_k*}(t,x)),\\
\ad u^{kj*}(t,x)=-(\sigma(t,s_{kj})\sigma'(t,s_{kj}))^{-1}B'(t,s_{kj})[x+(\lambda-z)
\lbar{H}(t,k)].\eea

In the following, we denote $n$th component of the optimal control
for this limit system as $u^{kj*}_n(t,x)$ Using such controls, we
construct \beq{con} u^\e(t,\al^\e(t),
x)=\sum^l_{k=1}\sum^{m_k}_{j=1}I_{\{\al^\e(t)=s_{kj}\}}
u^{kj*}(t,x).\eeq for the original problem. This control can also be
written as if $\al^\e(t)\in\mathcal{M}_k,u^\e(t,\al^\e(t),
x)=-(\sigma(t,\al^\e(t))\sigma'(t,\al^\e(t)))^{-1}B'(t,\al^\e(t))[x+(\lambda-z)
\lbar{H}(t,\lbar{\al}^\e(t))]$.  To proceed, we present the
following lemmas first.

\begin{lem}\label{re}
For a positive $T$ and any $k=1,2,\ldots,l,j=1,2,\ldots, m_k$,
\beq{ind1} \sup_{0\le t\le T}
E\l\int^t_0[I_{\{\al^\e(s)=s_{kj}\}}-\mu^k_j
I_{\{\lbar{\al}^\e(s)=k\}}]x^\e(s)r(s,s_{kj})ds\r^2\to 0 \text{ as }
\e \to 0.\eeq
\end{lem}

The proof is omitted for brevity.

\begin{lem}\label{le2}
For any $k=1,2,\ldots,l,j=1,2,\ldots, m_k$,
\beq{ind1a}
E(I_{\{\lbar{\al}^\e(s)=k\}}-I_{\{\lbar{\al}(s)=k\}})^2\to 0 \text{
as } \e \to 0.\eeq
\end{lem}

\para{Proof.} Similar to \cite[Theorem 7.30]{Yin&Zhang}, we can show
that
$(I_{\{\lbar{\al}^\e\cd=1\}},\ldots, I_{\{\lbar{\al}^\e\cd=l\}})$
converges weakly to
$(I_{\{\bar{\al}\cd=1\}},\ldots,I_{\{\bar{\al}\cd=l\}})$ in $(D[0,T]:
\rr^l)$
as $\e \to 0$. By means of Cram\'{e}r-Word's device, for each $i\in
\mathcal{M}$, $I_{\{\lbar{\al}^\e\cd=i\}}$ converges weakly to
$I_{\{\bar{\al}\cd=i\}}$. Then by virtue of the Skorohod
representation (with a slight abuse of notation), we may
assume $I_{\{\lbar{\al}^\e\cd=i\}}\to I_{\{\bar{\al}\cd=i\}}$ w.p.1.
without change of notation. Now by dominance convergence theorem, we
can conclude the proof. \qed

 \begin{thm}\label{1}
The control $u^\e(t)$ defined in \eqref{con} is nearly optimal in
that
 $\lim_{\e \to 0} |J^\e(\al,x,u^\e(\cdot))-v^\e(\al,x)|=0.$
\end{thm}

\para{Proof.}
Recall the definition of $\rho(t,s_{kj})$ in \eqref{r} and note that
the constructed control is given as
$u^\e(t,x,\al^\e(t))=-(\sigma(t,\al^\e(t))\sigma'(t,\al^\e(t)))^{-1}
B'(t,\al^\e(t))[x+(\lambda-z)
\overline{H}(t,\lbar{\al}^\e(t))]$. Then $x^\e(t)$ follows
\bea%
dx^\e(t)\ad=\sum^l_{k=1}\sum^{m_k}_{j=1}[r(t,s_{kj})x^\e(t)-\rho(t,s_{kj})x^\e(t)-\rho(t,s_{kj})(\lambda-z)\lbar{H}(t,k)]
I_{\{\al^\e(t)=s_{kj}\}}dt
\\
\aad \ +\sum^d_{i=1}\sqrt{\sum^l_{k=1}\sum^{m_k}_{j=1}(\sum^{d_1}_{n=1}u^\e_n(t,x^\e(t),\al^\e(t))\sigma_{ni}(t,\al^\e(t)))^2I_{\{\al^\e(t)=s_{kj}\}}}dw_i(t).\\
x^\e(0)\ad=\hat{x}.\eea %Here for the sake of simplicity of subsequent
The cost function $J^\e(\al,x,u^\e\cd)=E[x^\e(T)+\lambda-z]^2$. Let
$x^*(t)$ be the optimal trajectory of the limit problem. Recall the
definition of $f(\cdot)$ and $g(\cdot)$ in the \thmref{4.1}. Then
\bea dx^*(t)\ad=f(t, x^*(t), \lbar{\al}(t),
U^*(t))dt+\sum^d_{i=1}g_i(t,\lbar{\al}(t),U^*(t))dw_i(t), \
x^*(0)=\hat{x}.\eea
 Similar to the methods in
\cite[Theorem 9.8]{Yin&Zhang},  for all $\al \in
\mathcal{M}_k$, and $k=1,2,\ldots,l$,
$$\lim_{\e\to 0}v^\e(x,\al)=\lbar{v}(x,k).$$ Here $\lbar{v}(x,k)$ is the value function of the limit problem.
  For any
$\al \in \mathcal{M}_k, k=1,2,\ldots,l$,
\bea%
0\le
|J^\e(x,u^\e(\cdot),\al)-v^\e(x,\al)|=|J^\e(x,u^\e(\cdot),\al)-\lbar{v}(x,k)+\lbar{v}(x,k)-
v^\e(x,\al)|.\eea To establish the assertion, it suffices to show
that \bea |J^\e(x,u^\e(\cdot),\al)-\lbar{v}(x,k)|\to 0, \eea
\beq{val}%
\barray
|J^\e(x,u^\e(\cdot),\al)-\lbar{v}(x,\al)|\ad=|E[x^\e(T)+\lambda-z]^2-E[x^*(T)+\lambda-z]^2|\\
\ad=|Ex^\e(T)^2+2(\lambda-z)Ex^\e(T)-E{x^*}^2(T)-2(\lambda-z)Ex^*(T)|\\
\ad\leq CE^\frac{1}{2}[x^\e(T)-x^*(T)]^2 \earray\eeq  for some
constant $C$. Here, H\"{o}lder inequality and finite second moment
of $x^\e(T)$ and $x^*(T)$ are used. Note that we can write
$E(x^\e(T)-x^*(T))^2$ as follows:
\beq {xxxx}%
\barray%
\ad E(x^\e(T)-x^*(T))^2 \\
\aad \le
K\sum^l_{k=1}\sum^{m_k}_{j=1}E(\int^T_0[r(s,s_{kj})x^\e(s)(I_{\{\al^\e(s)=s_{kj}\}}-\mu^k_jI_{\{
\lbar{\al}^\e(s)=k\}})]ds)^2\\
\aad\ \
+K\sum^l_{k=1}\sum^{m_k}_{j=1}E(\int^T_0[\mu^k_jr(s,s_{kj})(x^\e(s)-x^*(s))I_{\{
\lbar{\al}^\e(s)=k\}}]ds)^2\\
\aad\ \
+K\sum^l_{k=1}\sum^{m_k}_{j=1}E(\int^T_0[\mu^k_jr(s,s_{kj})x^*(s)(I_{\{\lbar{\al}^\e(s)=k\}}-I_{\{\lbar{\al}(s)=k\}})]ds)^2\\
\aad\ \
-K\sum^l_{k=1}\sum^{m_k}_{j=1}E(\int^T_0\rho(s,s_{kj})x^\e(s)(I_{\{\al^\e(s)=s_{kj}\}}-\mu^k_jI_{\{\lbar{\al}^\e(s)=k\}})ds)^2\\
\aad\ \
+K\sum^l_{k=1}\sum^{m_k}_{j=1}E(\int^T_0[\mu^k_j\rho(s,s_{kj})(x^\e(s)-x^*(s))
I_{\{\lbar{\al}^\e(s)=k\}}]ds)^2\\
\aad\ \ +K\sum^l_{k=1}\sum^{m_k}_{j=1}
E(\int^T_0\mu^k_j\rho(s,s_{kj})x^*(s)
(I_{\{\lbar{\al}^\e(s)=k\}}-I_{\{\lbar{\al}(s)=k\}})ds)^2\\
\aad\ \
-K\sum^l_{k=1}\sum^{m_k}_{j=1}E(\int^T_0[\rho(s,s_{kj})(\lambda-z)\lbar{H}(s,k)(I_{\{\al^\e(s)=k\}}
-\mu^k_jI_{\{\lbar{\al}^\e(s)=k\}})]ds)^2\\
\aad\ \
+K\sum^l_{k=1}\sum^{m_k}_{j=1}E(\int^T_0[\rho(s,s_{kj})(\lambda-z)\lbar{H}(s,k)\mu^k_j(I_{\{\lbar{\al}^\e(s)=k\}}
-I_{\{\lbar{\al}(s)=k\}})]ds)^2 +D, \earray \eeq%
where \bea D=\ad
 KE\Bigg[\int^T_0\sum^d_{i=1}\Big[\sqrt{\sum^l_{k=1}\sum^{m_k}_{j=1}(\sum^{d_1}_{n=1}u^\e_n(s,x^\e(s),\al^\e(s))
 \sigma_{ni}(s,\al^\e(s)))^2I_{\{\al^\e(s)=s_{kj}\}}}\\
\aad \qquad \qquad
-\sqrt{\sum^l_{k=1}\sum^{m_k}_{j=1}{\mu^k_j}(\sum^{d_1}_{n=1}u^{kj*}_n(s,x^*(s),\lbar{\al}(s))\sigma_{ni}(s,\al^\e(s)))^2
I_{\{\lbar{\al}(s)=k\}}}\Big]dw_i(s)\Bigg]^2.\eea First, we use
\lemref{re}, \lemref{le2}, and H\"{o}lder inequality repeatedly to
handel the drift part. For the diffusion part, realizing that
\bea%
D\le
%\ad K E
%\int^T_0\sum^d_{i=1}\Big[\sqrt{\sum^l_{k=1}\sum^{m_k}_{j=1}(\sum^{d_1}_{n=1}u^\e_n(s,x^\e(s),\al^\e(s))
% \sigma_{ni}(s,\al^\e(s)))^2I_{\{\al^\e(s)=s_{kj}\}}}\\
%\aad \quad
%-\sqrt{\sum^l_{k=1}\sum^{m_k}_{j=1}{\mu^k_j}(\sum^{d_1}_{n=1}u^{kj*}_n(s,x^*(s),\lbar{\al}(s))
%\sigma_{ni}(s,\al^\e(s)))^2I_{\{\lbar{\al}(s)=k\}}}\Big]^2ds\\
%\leq
\ad KE\int^T_0\sum^{d}_{i=1}\Bigg[\sqrt{\sum^l_{k=1}
\sum^{m_k}_{j=1}(\sum^{d_1}_{n=1}u^\e_n(s,x^\e(s),\al^\e(s))
 \sigma_{ni}(s,\al^\e(s)))^2 [I_{\{\al^\e(s)=s_{kj}\}}-
 \mu^k_jI_{\{\lbar{\al}^\e(s)=k\}}]}\\
%\aad \quad+\sqrt{\sum^l_{k=1}
%\sum^{m_k}_{j=1}(\sum^{d_1}_{n=1}u^\e_n(s,x^\e(s),\al^\e(s))
% \sigma_{ni}(s,\al^\e(s)))^2\mu^k_jI_{\{\lbar{\al}^\e(s)=k\}}}\\
%\aad \quad -\sqrt{\sum^l_{k=1}\sum^{m_k}_{j=1}
%{\mu^k_j}(\sum^{d_1}_{n=1}u^{kj*}_n(s,x^*(s),\lbar{\al}(s))
% \sigma_{ni}(s,\al^\e(s)))^2I_{\{\lbar{\al}^\e(s)=k\}}}\\
\aad \qquad +\sqrt{\sum^l_{k=1}\sum^{m_k}_{j=1}
{\mu^k_j}(\sum^{d_1}_{n=1}u^{kj*}_n(s,x^*(s),\lbar{\al}(s))
 \sigma_{ni}(s,\al^\e(s)))^2
[I_{\{\lbar{\al}^\e(s)=k\}}-I_{\{\lbar{\al}(s)=k\}}}]\\
\aad \qquad+\(x^\e(s)-x^*(s)\)\Bigg]^2ds .\eea

Here, we plugged in the control constructed in \eqref{con} for the
last term above and utilized the non-degeneracy assumption mentioned
in the previous section. %Note that our control is of  feedback form
%and  the following fact that \beq{y} \barray\disp
%E\int^T_0(\sqrt{I_{\{\al^\e(s)=s_{kj}\}}-\mu^k_jI_{\{\lbar{\al}(s)=k\}}})^2ds
%\le K
%E|\int^T_0(I_{\{\al^\e(s)=s_{kj}\}}-\mu^k_jI_{\{\lbar{\al}(s)=k\}})|^2ds.
%\earray \eeq
Then we can use property of stochastic integral, dominance
convergence theorem, similar techniques involved in dealing with the
drift part and the finite second moment of $x^\e(\cdot)$ and
$x^*(\cdot)$ to proceed with the diffusion part. Finally, after
detailed calculation, we have $E(x^\e(T)-x^*(T))^2 \le
o(\e)+K\int^T_0E(x^\e(s)-x^*(s))^2ds.$ Now with the help of
Gronwall's inequality, we obtain $E(x^\e(T)-x^*(T))^2\to 0 \text{
as }\e\to 0$.  The proof is thus concluded. \qed

\begin{rem}\label{Rem1}{\rm
Note that efficient frontier, efficient portfolio, and minimum
variance for the limit
system can be obtained similar to
\cite[Therome 5.1-5.3]{Y&Z}.
We
provide the discussion
 below and omit the detailed proofs.
 The discussions below also carry over to the case in which the
 Markov chain has transient states to be studied in the next section.
\begin{itemize}
  \item[(1)] {\rm If \eqref{fe} holds, we have%
 \beq{ef} \barray
\lbar{P}(0,\al)\lbar{H}(0,\al)^2+\theta-1<0.\earray\eeq and the
efficient control corresponding to $z$ is%
 \beq{efu}\barray
u^{kj*}(t,x)=-(\sigma(t,s_{kj})\sigma'
(t,s_{kj}))^{-1}B'(t,s_{kj})[x+(\lambda^*-z)
\lbar{H}(t,k)]\earray \eeq in which%
 \beq{lam}\barray
 \disp\lambda^*-z=\frac{z-\lbar{P}
 (0,\al)\lbar{H}(0,\al)\hat{x}}{\lbar{P}
 (0,\al)\lbar{H}(0,\al)^2+\theta-1}.\earray\eeq We
 can further show that among all the
 flow of the network system
 satisfying that the expected terminal
 flow value is $z$, the optimal variance of $x(T)$ is%
\beq{var}\barray
E(x^*(T)-z)^2\ad=\frac{\lbar{P}(0,\al)
\lbar{H}(0,\al)^2+\theta}{1-\theta-\lbar{P}(0,\al)\lbar{H}(0,\al)^2}
[z-\frac{\lbar{P}(0,\al)\lbar{H}(0,\al)}
{\lbar{P}(0,\al)\lbar{H}(0,\al)^2+\theta}\hat{x}]^2\\
\aad+\frac{\lbar{P}(0,\al)\theta}
{\lbar{P}(0,\al)\lbar{H}(0,\al)^2+\theta}\hat{x}^2.\earray\eeq
Therefore, The minimum terminal variance is%
 \beq{minva}\barray
\disp
E(x^*(T)-z)^2=\frac{\lbar{P}(0,\al)\theta}{\lbar{P}(0,\al)\lbar{H}(0,\al)^2+\theta}\hat{x}^2
\ge 0\earray\eeq with the minimum expected terminal flow of the
system%
\beq{mex}\barray \disp
z_{\min}=\frac{\lbar{P}(0,\al)\lbar{H}(0,\al)}{\lbar{P}(0,\al)\lbar{H}(0,\al)^2+\theta}\hat{x}.\earray\eeq
and corresponding Lagrange multiplier $\lambda^*_{\min}=0$.}
\item[(2)] {\rm
Assume that an efficient portfolio $u^*_1(t)$ is given by \eqref{efu}
corresponding to $z=z_1>z_{\min}$. Then a control $u^*(t)$ is
efficient
if and only if there is a $\pi\ge 0$ such that%
 \beq{eu}\barray
 u^*(t)=(1-\pi)u^*_{\min}(t)+\pi u^*_1(t).\earray\eeq
 where $t\in[0,T]$ and%
  \beq{umin}\barray
u^*_{\min}(t)=-(\sigma(t,s_{kj})\sigma'(t,s_{kj}))^{-1}
B'(t,s_{kj})[x-z_{\min}
\lbar{H}(t,k)].\earray \eeq}

Assertion (2) is known as ``mutual fund theorem" in the financial market
problems. In platoon control problems,
this result offers a practical way of selecting the optimal flow
controls so that the total platoon length can
be as close to the designated value in the sense
that the variance of the platoon length
is minimized. Similarly, in platoon
communication resource allocation problems,
this strategy is optimal in the sense that the designated
total throughput for the platoon communication
network is most efficiently used.

\end{itemize}
}\end{rem}

\subsection{Inclusion of Transient States}
In this section, we consider the case in which the Markov chain has
transient states. We assume  $\wdt{Q}= \left( \begin{array}{l}
\wdt Q_1 \;\,0 \\
 \wdt Q_0 \;\wdt Q_*  \\
 \end{array} \right)$
 where
 $\wdt{Q}_1=\diag\{\wdt{Q}^1,\wdt{Q}^2,\ldots,\wdt{Q}^l\}$,
 $\wdt{Q}_0=(\wdt{Q}^1_*,\ldots,\wdt{Q}^l_*).$ For
 each $k=1,2,\ldots,l$, $\wdt{Q}^k$ is a generator with dimension
 $m_k\times m_k$, $\wdt{Q}_*\in \mathbb{R}^{m_*\times m_*},$ $\wdt{Q}^k_* \in \mathbb{R}^{m_*\times
 m_k}$ and
 $m_1+m_2+\cdots+
 m_*=m$. The state space of the underlying Markov chain is given by
 $\mathcal{M}=\mathcal{M}_1\cup\mathcal{M}_2\cup
 \cdots\cup\mathcal{M}_*=\{s_{11},\ldots,s_{1m_1},\ldots,s_{l1}\ldots,s_{lm_l},s_{*1},\ldots,s_{*m_*}
 \}$, where $\mathcal{M}_*=\{s_{*1},s_{*2},\ldots,s_{*m_*}\}$
 consists of the transient states. Suppose for $k=1,2,\ldots,l$,
 $\wdt{Q}^k$ are irreducible, and $\wdt{Q}_*$ is Hurwitz, i.e., all of its
 eigenvalues have negative real parts.
 Let
$ \wdh{Q}=
\left( \begin{array}{l}
 \wdh Q^{11} \;\wdh Q^{12}  \\
 \wdh Q^{21} \;\wdh Q^{22}  \\
 \end{array} \right)$
where $\wdh{Q}^{11}\in \mathbb{R}^{(m-m_*)\times(m-m*)}$,
$\wdh{Q}^{12}\in \mathbb{R}^{(m-m_*)\times m_*}$, $\wdh{Q}^{21}\in
\mathbb{R}^{m_*\times (m-m_*)}$, and $\wdh{Q}^{22}\in
\mathbb{R}^{m_*\times m_*}$. We define%
\bea
\lbar{Q}_*=\diag(\mu^1,\ldots,\mu^l)(\wdh{Q}^{11}\tilde{\mathbbold{1}}+\wdh{Q}^{12}(a_{m_1},a_{m_2},\ldots,a_{m_l})).
\eea with $\tilde{\mathbbold{1}}=\diag(\mathbbold{1}_{m_1},\ldots,
\mathbbold{1}_{m_l})$, $\mathbbold{1}_{m_j}=(1,\ldots, 1)'\in
\mathbb{R}^{m_j \times 1}$ and, for $k=1,\ldots,l$,
\bea%
a_{m_k}=(a_{{m_k},1},\ldots,a_{{m_k},m_*})'=-\wdt{Q}^{-1}_*\wdt{Q}^k_*\mathbbold{1}_{m_k}.\eea
Let $\xi$ be a random variable uniformly distributed on $[0,1]$ that
is independent of $\al^\e\cd$. For each $j=1,2,\ldots,m_*$, define
an integer-valued random variable $\xi_j$ by \bea%
\xi_j=I_{\{0 \le\xi\le a_{{m_1},j}\}}+2I_{\{a_{{m_1},j} <\xi \le
a_{{m_1},j}+ a_{{m_2},j} \}}+\cdots+lI_{\{ a_{{m_1},j}+\cdots+
a_{{m_{l-1}},j}<\xi\le 1\}}.\eea Now define the aggregated process
$\lbar{\al}^\e\cd$ by \bea%
 \lbar{\al}^\e(t)= \left\{ \begin{array}{l}
 k,\;\text{ if }\al ^\e(t)\in {\cal M}_k  \\
 \xi _j ,\text{ if }\al ^\e (t) = s_{*j.}  \\
 \end{array} \right.
\eea Note the state space of $\lbar{\al}^\e(t)$ is $\overline{\cal
M}=\{1,2,\ldots,l\}$ and $\lbar{\al}^\e\cd\in D[0,T]$. In addition,
$$P(\lbar{\al}^\e(t)=i|\al^\e(t)=s_{*j})=a_{{m_i},j}.$$ Then
according to \cite[Theorem 4.2]{GZB}, $\lbar{\al}^\e(\cdot)
$ converges weakly to
$\lbar{\al}(\cdot)$ such that $\lbar{\al}(\cdot)\in
\{1,2,\ldots,l\}$ is a Markov chain generated by $\lbar{Q}_*$.

\begin{thm}
As $\e \to 0$, we have $P^\e(s,s_{kj})\to \overline{P}(s,k)$ and
$H^\e(s,s_{kj})\to \overline{H}(s,k)$, for $k=1,2,\ldots, l$,
$j=1,2,\ldots,m_k$, $P^\e(s,s_{*j})\to \overline{P}_*(s,j)$ and
$H^\e(s,s_{*j})\to \overline{H}_*(s,j)$, for $j=1,2,\ldots, m_*$
uniformly on $[0,T]$ where \bea
\overline{P}_*(s,j)=a_{{m_1},j}\overline{P}(s,1)+\cdots+a_{{m_l},j}\overline{P}(s,l)\eea
\bea
\overline{H}_*(s,j)=a_{{m_1},j}\overline{H}(s,1)+\cdots+a_{{m_l},j}\overline{H}(s,l)\eea
and $\overline{P}(s,k)$ and $\overline{H}(s,k)$ are the unique
solutions to the following equations. For $k=1,2,\ldots,l$,
\beq{P2T} \barray
\dot{\overline{P}}(t,k)=\ad(\overline{\rho}(t,k)-2\overline{r}(t,k))\overline{P}(t,k)-\overline{Q}_*\overline{P}(t,\cdot)(k)\\
\overline{P}(T,k)=\ad 1\earray \eeq and \beq{H2a}\barray
\dot{\overline{H}}(t,k)=\ad\overline{r}(t,k)\overline{H}(t,k)
-\frac{1}{\overline{P}(t,k)}\overline{Q}_*\overline{P}
(t,\cdot)\overline{H}(t,\cdot)(k)
+\frac{\overline{H}(t,k)}{\overline{P}(t,k)}\overline{Q}_*\overline{P}(t,\cdot)(k)\\
\overline{H}(T,k)=\ad 1. \earray\eeq
\end{thm}

The convergence of $P^\e$ and $H^\e$ leads to $v^\e(t,s_{kj},x) \to
\lbar{v}(t,k,x)$, for $k=1,2,\ldots,l,j=1,2,\ldots,m_k$,
$v^\e(t,s_{*j},x)\to v_*(t,j,x)$ for $j=1,2,\ldots,m_*$, where\bea
v_*(t,j,x)=a_{{m_1},j}\lbar{v}(t,1,x)+\cdots+a_{{m_l},j}\lbar{v}(t,l,x)\eea
and $\lbar{v}(t,k,x)=\lbar{P}(t,k)(x+(\lambda-z)\lbar{H}(t,k))^2$.
The control set for the limit problem is the same as that for
recurrent case and is given by \bea
\mathcal{U}=\{U=(U^1,U^2,\ldots,U^l):
U^k=(u^{k1},u^{k2},\ldots,u^{km_k}),u^{kj}\in
\mathbb{R}^{d_1}\}.\eea Then the corresponding
limit problem is \bea dx(t)=f(x(t),\lbar{\al}(t),
U(t))dt+\sum^d_{i=1}g_i(t,\lbar{\al}(t),U(t))dw_i(t).\eea where
$\overline{\al}(\cdot)\in \{1,2,\ldots,l\}$ is a Markov chain
generated by $\lbar{Q}_*$ with $\lbar{\al}(0)=\al$. The optimal
control for this limit problem is \bea
U^*(t)=(U^{1*}(t,x),U^{2*}(t,x),\ldots,U^{l*}(t,x)).\eea with
\[U^{k*}(t,x)=(u^{k1*}(t,x),u^{k2*}(t,x),\ldots,u^{km_k*}(t,x))\] and
\[u^{kj*}(t,x)=-(\sigma(t,s_{kj})\sigma'(t,s_{kj}))^{-1}B'(t,s_{kj})[x+(\lambda-z)
\lbar{H}(t,k)].\] Using such controls, we construct \beq{con2}
u^\e(t,\al^\e(t),
x)=\sum^l_{k=1}\sum^{m_k}_{j=1}I_{\{\al^\e(t)=s_{kj}\}}
u^{kj*}(t,x)+\sum^{m_*}_{j=1}I_{\{\al^\e(t)=s_{*j}\}}u^{*j*}(t,x).\eeq
for the original problem where
$u^{*j*}(s,x)=-(\sigma(t,s_{*j})\sigma'(t,s_{*j}))^{-1}B'(t,s_{*j})[x+(\lambda-z)
\lbar{H}_*(s,j)]$.

\para{Proof.}
Following the proof of \thmref{4.1} to \eqref{ct},  we have for
$s\in[0,T],$ $\wdt{Q}^kP^0(s,\cdot)(s_{kj})=0$ for $k=1,2,\ldots,l$,
$j=1,2,\ldots,m_k$%
\bea%
\ad(\wdt{Q}^1_*,\ldots,\wdt{Q}^l_*,\wdt{Q}_*)(P^0(s,s_{11}),\ldots,P^0(s,s_{1m_1}),
\ldots,\\
\aad \quad \times P^0(s,s_{l1}),\ldots,P^0(s,s_{lm_l}),P^0(s,s_{*1}),\ldots,P^0(s,s_{*m*}))'=0.%
\eea%
The irreducibility of $\wdt{Q}^k$ for any $k$ implies \bea
(P^0(s,s_{k1}),\ldots,P^0(s,s_{k_{m_k}}))'=P^0(s,k)\mathbbold{1}_{m_k}.\eea
Let $P_*(s)=(P^0(s,s_{*1}),\ldots, P^0(s,s_{*m*}))'$, we have \bea
\wdt{Q}^1_*\mathbbold{1}_{m_1}P^0(s,1)+\cdots+\wdt{Q}^l_*\mathbbold{1}_{m_l}P^0(s,l)+\wdt{Q}_*P_*(s)=0.\eea
Here, \bea
P_*(s)=-\wdt{Q}^{-1}_*(\wdt{Q}^1_*\mathbbold{1}_{m_1}P^0(s,1)+\cdots+\wdt{Q}^l_*\mathbbold{1}_{m_l}P^0(s,l))
=a_{m_1}P^0(s,1)+\cdots+a_{m_l}P^0(s,l).\eea Then $P_*(s)\in
\rr^{m_*}$ and its $j$th component is $P_*(s,j)$.
 The rest of the proof is
similar to that of \thmref{4.1}, except replacing $\lbar{Q}$ by
$\lbar{Q}_*$.

\begin{thm}\label{thm:2}
The control $u^\e(t)$ defined in \eqref{con2} is nearly optimal in
that
$\lim_{\e \to 0} |J^\e(\al,x,u^\e\cd)-v^\e(\al,x)|=0.$
\end{thm}

\para{Proof.}
The proof is similar to that of Theorem \thmref{1} with the
use of the estimate
$ E|\int^t_0I_{\{\al^\e(s)=s_{*j}\}}ds|^2
\to 0$  as  $\e \to 0$  from \cite[Theorem 3.1]{GZB}. \qed

 \section{Further Remarks}\label{sec:rem}
 This work focused on the near-optimal controls for non-definite
 control problems. Our primary motivation stems from networked systems.
 Our approach provides a
 systematic approach to reduce the complexity of the underlying
 system. In lieu of treating the large dimensional systems directly, we
 solve a set of limit Riccati equations that have much
 smaller dimensions.
Using the limit problems as a guide to design controls for the
original systems leads to near-optimal controls of the original
systems. Although the paper is devoted to platoon controls, the
results can be readily applied to the portfolio optimization in
financial engineering. Future research efforts
 can be directed to the study of non-definite
 control problems in the hybrid systems,
 in which the Markov chain $\alpha(t)$
 is a hidden process and
 Wonham filter will be involved.  More thoughts and further considerations
are needed.

\end{document}